\documentclass[11pt]{amsart}
\usepackage{hyperref}
\usepackage{amssymb,amsfonts}
\usepackage{hyperref}
\usepackage{amssymb,textcomp}

\usepackage{enumerate}

\usepackage[utf8]{inputenc}
\usepackage[T1]{fontenc}

\usepackage[mathscr]{eucal}

  \newcommand{\B}{\mathcal{B}}
 
\newcommand{\cX}{\mathcal{X}}

\newcommand{\C}{\mathbb{C}}
\newcommand{\T}{\mathbb{T}}
\newcommand{\R}{\mathbb{R}}

\newcommand{\N}{\mathbb{N}}
\newcommand{\bbZ}{\mathbb{Z}}
\newcommand{\Z}{\mathbb{Z}}

\newcommand{\norm}[1]{\left\Vert #1\right\Vert}

\theoremstyle{plain}
\newtheorem{theorem}{Theorem}[section]

\newtheorem{problem}{Problem}

\newtheorem*{theorem*}{Theorem}

\theoremstyle{definition}
\newtheorem{definition}[theorem]{Definition}

\theoremstyle{remark}

\begin{document}
\title{Ergodic theory: Recurrence}

\author{Nikos Frantzikinakis}
\address[Nikos Frantzikinakis]{University of Crete, Department of Mathematics, Voutes University Campus, Heraklion 71003, Greece} \email{frantzikinakis@gmail.com}
\author{Randall McCutcheon}
\address[Randall McCutcheon]{University of Memphis, Department of Mathematics, Memphis, TN 38152-3240, USA} \email{rmcctchn@memphis.edu}

\thanks{The first author was partially supported
 by NSF grant DMS-0701027.}

\maketitle
\setcounter{tocdepth}{1}
\tableofcontents


\noindent {\bf Almost every, essentially}: Given a Lebesgue
measure space $(X,\mathcal{B},\mu)$, a property $P(x)$ predicated
of elements of $X$ is said to hold for almost every $x\in X$, if
the set $X\setminus  \{x\colon P(x) \text { holds}\}$ has zero
measure. Two sets $A,B\in \mathcal{B}$ are essentially disjoint if
$\mu(A\cap B)=0$.

\noindent{\bf Conservative system}: Is an infinite measure
preserving system such that for no set $A\in\mathcal{B}$ with
positive measure are $A, T^{-1}A, T^{-2}A,\ldots$ pairwise
essentially disjoint.

\noindent{\bf $(c_n)$-conservative system}:
If $(c_n)_{n\in \N}$ is a
decreasing sequence of positive real numbers, a conservative
ergodic measure preserving transformation $T$ is
$(c_n)$-conservative if for some non-negative function $f\in
L^1(\mu)$, $\sum_{n=1}^\infty c_n f(T^nx)=\infty$ a.e.

\noindent{\bf Doubling map}: If  $\T$  is the interval $[0,1]$
with its endpoints identified and addition performed modulo $1$,
the (non-invertible) transformation $T\colon \T\to \T$, defined by
$Tx=2x$ mod 1, preserves Lebesgue measure, hence induces a measure
preserving system on $\T$.

\noindent{\bf Ergodic system}:  Is a measure preserving system
$(X, \mathcal{B},\mu, T)$ (finite or infinite) such that every
$A\in\mathcal{B}$  that is  $T$-invariant (i.e. $T^{-1}A=A$)
satisfies either $\mu(A)=0$ or $\mu(X\setminus A)=0$. (One can
check that the rotation $R_\alpha$ is ergodic if and only if
$\alpha$ is irrational, and that the doubling map is ergodic.)

\noindent{\bf Ergodic decomposition}: Every measure preserving
system $(X,\cX,\mu,T)$ can be expressed as an integral of ergodic
systems; for example, one can write $\mu=\int \mu_t\ d\lambda(t)$,
where $\lambda$ is a probability measure on $[0,1]$ and $\mu_t$ are
$T$-invariant probability measures on $(X,\cX)$ such that the
systems $(X,\cX,\mu_t,T)$ are ergodic for $t\in [0,1]$.

\noindent{\bf Ergodic theorem}:
States that if $(X,\mathcal{B},\mu,T)$ is a measure
preserving system and $f\in L^2(\mu)$, then
$\lim_{N\to\infty}\norm{\frac{1}{N}\sum_{n=1}^NT^nf-P_f}_{L^2(\mu)}=0$,
where $P_f$ denotes the orthogonal projection of the function $f$
onto the subspace $\{f \in L^2(\mu)\colon Tf=f\}$.

\noindent {\bf Hausdorff $a$-measure}:
Let $(X,\mathcal{B},\mu,T)$ be a measure preserving system endowed
with a $\mu$-compatible metric $d$. The Hausdorff
$a$-measure $\mathcal{H}_a(X)$ of $X$ is an outer measure
defined for all subsets of $X$ as follows: First, for $A\subset X$
and $\varepsilon>0$ let $\mathcal{H}_{a,\varepsilon}(A)=\inf
\{\sum_{i=1}^\infty r_i^a\}$, where the infimum is taken over all
countable coverings of $A$ by sets $U_i\subset X$ with diameter
$r_i<\varepsilon$. Then define
$\mathcal{H}_a(A)=\limsup_{\varepsilon\to 0} \mathcal{H}_{a,
\varepsilon}(A)$.

\noindent{\bf Infinite measure preserving system}: Same as measure
preserving system, but $\mu(X)=\infty$.

\noindent{\bf Invertible system}: Is a measure preserving system
$(X, \mathcal{B},\mu,T)$  (finite or infinite), with the property
that there exists $X_0\in X$, with $\mu(X\setminus X_0)=0$, and
such that the transformation $T\colon X_0\to X_0$ is bijective,
with $T^{-1}$ measurable.

\noindent{\bf Measure  preserving system}:
Is a quadruple $(X, \mathcal{B},\mu,T)$, where $X$ is a set,
$\mathcal{B}$ is  a $\sigma$-algebra  of subsets of $X$ (i.e.
$\mathcal{B}$ is closed under countable unions and
complementation), $\mu$ is a probability measure (i.e. a countably
additive function from $\mathcal{B}$ to $[0,1]$ with $\mu(X)=1$),
and $T\colon X \to  X$ is measurable (i.e. $T^{-1}A=\{ x\in X:
Tx\in A\}\in \mathcal{B}$ for $A\in\mathcal{B}$), and
$\mu$-preserving  (i.e. $\mu(T^{-1} A)=\mu(A)$). Moreover,
throughout the discussion we assume that the measure space
$(X,\mathcal{B},\mu)$ is Lebesgue (see Section $1.0$ of
\cite{Aa}).

\noindent {\bf $\mu$-compatible metric}: Is a separable metric
on $X$, where $(X,\mathcal{B},\mu)$ is a probability space,
having the property that open sets measurable.

\noindent{\bf Positive definite sequence}: Is a complex-valued
sequence $(a_n)_{n\in \Z}$ such that for any $n_1,\ldots , n_k\in
\Z$ and $z_1,\ldots ,z_k\in \C$, $\sum_{i,j=1}^k a_{n_i-n_j} z_i
\overline{z_j} \geq 0 $.

 \noindent{\bf Rotations on
$\T$}: If $\T$ is the interval $[0,1]$ with its endpoints
identified and addition performed modulo $1$, then for every
$\alpha \in \R$ the transformation $R_\alpha\colon \T\to \T$,
defined by $R_\alpha x=x+\alpha$, preserves Lebesgue measure on
$\T$ and hence induces a measure preserving system on $\T$.

\noindent{\bf Syndetic set}: Is a subset $E\subset \Z$ having
bounded gaps. If $G$ is a general discrete group, a set
$E\subset G$ is syndetic if $G=FE$ for some finite set $F\subset G$.

\noindent{\bf Upper density}: Is the number $\overline{d}(\Lambda)=
\limsup_{N\to\infty} {|\Lambda\cap\{-N,\ldots ,N\}|\over 2N+1}$, where
$\Lambda\subset\Z$ (assuming the limit to exist).
Alternatively for measurable $E\subset\R^m$,
$\overline{D}(E)=\limsup_{l(S)\to\infty}{ m(S\cap E)\over m(S)}$, where
$S$ ranges over all cubes in $\R^m$, and $l(S)$ denotes the length
of the shortest edge of $S$.

\noindent{\bf Notation}: The following notation will be used
throughout the article: $Tf=f\circ T$, $\{x\}=x-[x]$,
$\text{D-}\!\lim_{n\to\infty}(a_n)=a$ $\leftrightarrow$
$\overline{d} \big( \{ n: |a_n-a|>\varepsilon\} \big) =0$ for
every $\varepsilon>0$.

\section{Definition of the Subject and its Importance}
The basic principle that lies behind several recurrence
phenomena, is that the typical trajectory of a system with finite
volume comes back infinitely often to any neighborhood of its
initial point. This principle was first exploited by Poincar\'e
in his 1890 King Oscar prize-winning memoir that studied planetary
motion. Using the prototype of an ergodic theoretic argument, he
showed that in any system of point masses having fixed total
energy that restricts its dynamics to bounded subsets of its phase
space, the typical  state of motion (characterized by
configurations and velocities) must recur to an arbitrary degree
of approximation.

Among the recurrence principle's more spectacularly
counterintuitive ramifications is that isolated ideal gas systems
that do not lose energy will return arbitrarily closely to their
initial states, even when such a return entails a decrease in
entropy from equilibrium, in apparent contradiction to the second
law of thermodynamics. Such concerns, previously canvassed by
Poincar\'e himself, were more infamously expounded by  Zermelo
(\cite{Ze}) in 1896. Subsequent clarifications by Boltzmann,
Maxwell and others led to an improved understanding of the second
law's primarily statistical nature. (For an interesting
historical/philosophical discussion, see \cite{Skl}; also \cite{Be11}.
For a probabilistic analysis of the likelihood of observing second law
violations in small systems over short time intervals, see
\cite{SE}.)

 These discoveries had a profound impact in dynamics, and
 the theory of measure preserving transformations (ergodic theory)
 evolved from these developments. Since then,
 the Poincar\'e recurrence principle has been applied to a variety
 of different fields in mathematics, physics, and information
 theory. In this article we survey the impact
 it has had in ergodic theory, especially as pertains to the field
of {\it ergodic Ramsey theory.} (The heavy emphasis herein
on the latter reflects authorial interest, and is
not intended to transmit a proportionate image of the broader landscape
of research relating to recurrence in ergodic theory.)
Background information we assume in this article
can be found in the books \cite{EW11, Fu2,Gl03, HK17, Pe,
Wa}  (see also the article ``Measure Preserving Systems''
by K. Petersen in this volume). Related information can also be found on the survey articles \cite{Be1, Be2, Be06b, Fr16, Kra06, Kra07, Kra11}.

\section{Introduction}
In this section  we shall give several formulations  of the
Poincar\'e recurrence principle using the language of ergodic
theory. Roughly speaking, the principle states that in a finite (or
conservative) measure preserving system, every set of positive
measure (or almost every point) comes back
to itself infinitely many times under iteration.
Despite the profound importance of these results, their proofs are
extremely simple.
\begin{theorem}[Poincar\'e Recurrence for Sets]\label{PoincareI}
 Let
$(X,\mathcal{B},\mu, T)$ be a measure preserving
system  and $A\in \mathcal{B}$ with
$\mu(A)>0$. Then $\mu(A\cap T^{-n}A)>0$ for infinitely many $n\in
\N$.
\end{theorem}
\begin{proof}
Since $T$ is measure preserving, the sets
$A,T^{-1}A,T^{-2}A,\ldots$  have the same measure. These sets
cannot be pairwise essentially disjoint, since then the union of
finitely many of them  would have measure greater than $\mu(X)=1$.
Therefore, there exist $m,n\in \N$, with $n>m$, such that
$\mu(T^{-m}A\cap T^{-n}A)>0$. Again since
$T$ is measure preserving, we conclude that $\mu(A\cap T^{-k}A)>0$,
where $k=n-m>0$. Repeating this argument for the iterates
$A,T^{-m}A, T^{-2m}A,\ldots$, for all $m\in\N$, we easily deduce
that $\mu(A\cap T^{-n}A)>0$ for infinitely many $n\in\N$.
\end{proof}
We remark that the above argument actually shows that $\mu(A\cap
T^{-n}A)>0$ for some $n\leq [\frac{1}{\mu(A)}]+1$.
\begin{theorem}[Poincar\'e Recurrence for Points]\label{PoincareII} Let
$(X,\mathcal{B},\mu, T)$ be a measure preserving system and $A\in
\mathcal{B}$. Then for almost every  $x\in A$ we have that $T^nx
\in A$ for infinitely many $n\in \N$.
\end{theorem}
\begin{proof}
 Let $B$ be the set of $x\in A$
such that $T^nx\notin A$ for all $n\in \N$. Notice that
$B=A\setminus \bigcup_{n\in\N} T^{-n}A$; in particular, $B$ is measurable.
Since  the iterates
$B,T^{-1}B, T^{-2}B,\ldots$ are pairwise essentially disjoint, we
conclude (as in the proof of Theorem~\ref{PoincareI}) that
$\mu(B)=0$. This shows that for almost every  $x\in A$ we have
that $T^nx \in A$ for some $n\in \N$. Repeating this argument for
the transformation $T^m$ in place of $T$ for all $m\in\N$, we
easily deduce the advertised statement.
\end{proof}
Next we give a variation of Poincar\'e recurrence for  measure
preserving systems endowed with a compatible metric:
\begin{theorem}[Poincar\'e Recurrence for Metric Systems]\label{PoincareIII}
Let $(X,\mathcal{B},\mu, T)$ be a measure preserving system, and
suppose that $X$ is endowed with a
$\mu$-compatible metric. Then for almost every  $x\in X$ we have
$$\liminf_{n\to\infty} d(x,T^nx)=0.
$$
\end{theorem}

The proof of this result is similar to the proof of
Theorem~\ref{PoincareII} (see \cite{Fu2}, page 61). Applying  this
result to the doubling map $Tx=2x$ on $\T$,
we get that for almost every $x\in X$, every string of
zeros and ones in the dyadic expansion of $x$ occurs infinitely
often.

We remark that all three formulations of the Poincar\'e Recurrence
Theorem that we have given hold for conservative systems as well.
See, e.g., \cite{Aa} for details.


This article is structured as follows. In
Section~\ref{quantitative} we give a few quantitative versions of
the previously mentioned qualitative results. In
Sections~\ref{single} and \ref{S:multiple} we give several
refinements of the Poincar\'e recurrence theorem, by restricting
the scope of the return time $n$, and by considering multiple
intersections (for simplicity we focus on $\Z$-actions). In
Section~\ref{applications} we give various implications of the
recurrence results in combinatorics and number theory (see also
the article ``Ergodic Theory: Interactions with Combinatorics and
Number Theory'' by T. Ward in the present volume). Lastly, in
Section~\ref{problems} we give several open problems related to
the material presented in
Sections~\ref{single}-\ref{applications}.

\section{Quantitative Poincar\'e Recurrence}\label{quantitative}
\subsection{Early results} For applications it is desirable to have quantitative
versions of the results mentioned in the previous section.  For
example one would like to know  how large $\mu(A\cap T^{-n}A)$
can be made and for how many $n$.
\begin{theorem}[Khintchine \cite{Ki}]\label{Khin}
 Let $(X,\mathcal{B},\mu,T)$ be a measure preserving system
and $A\in \mathcal{B}$. Then for every
$\varepsilon>0$ we have  $\mu(A\cap T^{-n}A)>\mu(A)^2-\varepsilon$
for a set of $n\in\N$ that has bounded gaps.
\end{theorem}
By considering the doubling map $Tx=2x$  on $\T$ and  letting
$A={\bf 1}_{[0,1/2)}$, it is easy to check that the lower bound of
the previous result cannot be improved.
We also remark that it is not possible to estimate  the size of
the gap by a function of $\mu(A)$ alone. One can see this by
considering the rotations $R_{k}x=x+1/k$ for $k\in \N$, defined on
$\T$, and letting $A={\bf 1}_{[0,1/3]}$.

Concerning the second version of the Poincar\'e recurrence
theorem, it is natural to ask whether for almost every  $x\in X$
the set of return times $S_{x}=\{n\in\N\colon T^nx\in A\}$  has
bounded gaps. This is not the case, as one can see by considering
the doubling map $Tx=2x$ on $\T$ with the Lebesgue measure, and
letting $A={\bf 1}_{[0,1/2)}$. Since Lebesgue almost every  $x\in
\T$ contains arbitrarily large blocks of ones  in its dyadic
expansion, the set $S_x$ has unbounded gaps. Nevertheless, as an
easy consequence of the Birkhoff ergodic theorem (\cite{Bi}), one
has the following:
\begin{theorem} Let $(X,\mathcal{B},\mu,T)$ be a measure preserving system and
$A\in \mathcal{B}$ with $\mu(A)>0$. Then for almost every  $x\in
X$ the set $S_x=\{n\in\N\colon T^nx\in A\}$ has well defined
density and $\int d(S_x)\; d\mu(x) =\mu(A)$. Furthermore, for ergodic
measure preserving systems we have $d(S_x)= \mu(A)$ a.e.
\end{theorem}
Another question that arises naturally is, given a set $A$ with
positive measure and an $x\in A$, how long should one wait till
some iterate $T^nx$ of $x$ hits $A$? By considering an irrational
rotation  $R_\alpha$ on $\T$, where  $\alpha$ is very near to, but
not less than, ${1\over 100}$, and
letting $A=1_{[0,1/2]}$, one can see that the first return time
is a member of the set $\{1,50,51\}$. So it may come as a surprise that
the average first
return time does not depend on the system (as long as it is
ergodic), but only on the measure of the set $A$.
\begin{theorem}[Kac \cite{Ka}]
 Let $(X,\mathcal{B},\mu,T)$ be an ergodic measure preserving system
and $A\in \mathcal{B}$ with
 $\mu(A)>0$. For $x\in X $ define $R_A(x)=\min\{n\in\N\colon T^nx\in
 A\}$. Then for $x\in A$ the expected value of $R_A(x)$ is $1/\mu(A)$, i.e.
 $\int_A R_A(x) \ d\mu=1$.
\end{theorem}

\subsection{More Recent Results}
As we mentioned in the previous section, if the space $X$ is
endowed with a $\mu$-compatible  metric $d$, then for almost every
$x\in X$ we have that $\liminf_{n\to\infty} d(x,T^nx)=0$. A
natural question is, how much iteration is needed to come back
within a small distance of a given typical point? Under some
additional hypothesis on the metric $d$ we have the following
answer:
\begin{theorem}[Boshernitzan~\cite{Bos}]\label{T:Bos}
Let $(X,\mathcal{B},\mu,T)$ be a measure preserving system endowed
with a $\mu$-compatible metric $d$. Assume that the Hausdorff
$a$-measure $\mathcal{H}_a(X)$ of $X$ is
$\sigma$-finite (i.e., $X$ is a countable union of sets $X_i$ with
$\mathcal{H}_a(X_i)<\infty$). Then for almost every $x\in X$,
$$
\liminf_{n\to\infty} \big\{n^{\frac{1}{a}}\cdot d(x,
T^nx)\big\}<\infty.
$$
Furthermore, if $\mathcal{H}_a(X)=0$, then for almost every
$x\in X$,
$$
\liminf_{n\to\infty} \big\{n^{\frac{1}{a}}\cdot d(x,
T^nx)\big\}=0.
$$
\end{theorem}
One can see from rotations by ``badly approximable'' vectors
$\alpha\in \T^k$ that the exponent $1/k$ in the previous theorem
cannot be improved. Several applications of Theorem~\ref{T:Bos} to
billiard flows, dyadic transformations, symbolic flows and
interval exchange transformations are given in \cite{Bos}. For a
related result dealing with mean values of the limits  in
Theorem~\ref{T:Bos} see \cite{Sk}.

An interesting connection between rates of recurrence and entropy
of an ergodic measure preserving system  was established by
Ornstein and Weiss~(\cite{OW}), following earlier work of Wyner
and Ziv~(\cite{WZ}):
\begin{theorem}[Ornstein \& Weiss~\cite{OW}]
Let $(X,\mathcal{B},\mu,T)$ be an ergodic measure preserving
system and $\mathcal{P}$
be  a finite partition of  $X$. Let $P_n(x)$ be the
element of the partition
$\bigvee_{i=0}^{n-1}T^{-i}\mathcal{P}=\{ \bigcap_{i=0}^{n-1}T^{-i}
P^{(i)}: P^{(i)}\in \mathcal{P}, 0\leq i <n\}$ that contains
$x$. Then for almost every $x\in X$, the first return time $R_n(x)$  of
$x$ to $P_n(x)$ is
asymptotically equivalent to $e^{h(T,\mathcal{P})n}$, where
$h(T,\mathcal{P})$ denotes the entropy of the system with
respect to the partition $\mathcal{P}$. More precisely,
$$
\lim_{n\to\infty} \frac{\log{R_n(x)}}{n}=h(T,\mathcal{P}).
$$
\end{theorem}
An extension of the above result to some classes of infinite
measure preserving systems was given in \cite{Ga}.

Another connection of recurrence rates, this time with the local
dimension of an invariant measure, is given by the next result:
\begin{theorem}[Barreira~\cite{Ba1}]\label{T:Ba}
Let $(X,\mathcal{B},\mu,T)$ be an ergodic measure preserving
system. Define the upper and lower recurrence rates
$$
\underline{R}(x)=\liminf_{r\to 0 }
\frac{\log{\tau_r(x)}}{-\log{r}} \quad \text{and}\quad
\overline{R}(x)=\limsup_{r\to 0 }
\frac{\log{\tau_r(x)}}{-\log{r}},
$$
where $\tau_r(x)$ is the first return time of $T^kx$ in $B(x,r)$,
and the upper and lower pointwise dimensions
$$
\underline{d}_\mu(x)=\liminf_{r\to 0}
\frac{\log{\mu(B(x,r))}}{\log{r}} \quad  \text{and} \quad
\overline{d}_\mu(x)=\limsup_{r\to 0}
\frac{\log{\mu(B(x,r))}}{\log{r}}.
$$
Then  for almost every $x\in X$, we have
$$
\underline{R}(x)\leq \underline{d}_\mu(x) \quad \text{and} \quad
\overline{R}(x)\leq \overline{d}_\mu(x).
$$
\end{theorem}
Roughly speaking, this theorem asserts that for typical $x\in X$
and  for small $r$, the first return time of $x$ in $B(x,r)$ is at
most $r^{-d_\mu(x)}$. Since $\underline{d}_\mu(x)\leq
\mathcal{H}_a(X)$ for almost every $x\in X$, we can conclude the
first part of Theorem~\ref{T:Bos} from Theorem~\ref{T:Ba}. For
related results the interested reader should consult the survey
\cite{Ba2} and the bibliography therein.

We also remark that the  previous results and related concepts
have been applied to estimate the dimension of certain strange
attractors (see \cite{HJ} and the references therein) and the
entropy of certain Gibbsian systems~\cite{CU}.

We end this section with a result that connects ``wandering
rates'' of sets in infinite measure preserving systems with their
``recurrence rates''. The next  theorem follows easily from a
result about lower bounds on ergodic averages for measure
preserving systems due to Leibman (\cite{Le0}); a weaker form for
conservative, ergodic systems can be found in Aaronson
(\cite{Aa0}).

\begin{theorem}
Let $(X,\mathcal{B},\mu,T)$ be an infinite measure preserving
system, and $A\in\mathcal{B}$ with $\mu(A)<\infty$. Then for all
$N\in \N$,
$$
\Big( \frac{ \mu(\bigcup_{n=0}^{N-1} T^{-n}A)}{N}\cdot
\sum_{n=0}^{N-1} \mu(A\cap T^{-n}A) \Big) \geq \frac{1}{2}\cdot
(\mu(A))^2.
$$

 \end{theorem}
\section{Subsequence Recurrence}\label{single}
In this section we discuss  what  restrictions we can impose on
the set of return times in the various versions of the Poincar\'e
recurrence theorem. We start with:
\begin{definition}
Let $R\subset \Z$. Then $R$ is a set of:

$(a)$ {\em Recurrence}  if for any invertible measure preserving
system $(X,\mathcal{B},\mu,T)$, and $A\in \mathcal{B}$ with
$\mu(A)>0$, there is some nonzero $n\in R$ such that $\mu(A\cap
T^{-n}A)>0$.

$(b)$ {\em Topological recurrence} if for every compact metric
space $(X, d)$, continuous transformation $T\colon X\rightarrow X$
and every $\varepsilon>0$, there are $x\in X$ and nonzero $n\in
R$ such that $d(x,T^nx)<\varepsilon$.
\end{definition}
It is easy to check that the existence of a single $n\in R$
satisfying the previous recurrence  conditions actually guarantees
the existence of infinitely many $n\in R$ satisfying the same
conditions. Moreover, if $R$ is a set of recurrence then one can
see from existence of some $T$-invariant measure $\mu$
that $R$ is also a set of topological recurrence. A (complicated)
example showing that the converse is not true was given by
Kriz~(\cite{Kr}).

Before giving a list of examples of sets of (topological)
recurrence, we discuss some necessary conditions: A set of
topological recurrence must contain infinitely many multiples of
every positive integer, as one can see by considering rotations on
$\Z_d$, $d\in \N$ . Hence, the sets $\{2n+1,n\in \N\}$,
$\{n^2+1,n\in \N\}$, $\{p+2, p \text{ prime}\}$ are not good for
(topological) recurrence. If $(s_n)_{n\in \N}$ is a lacunary
sequence (meaning $\liminf_{n\to\infty}(s_{n+1}/s_n)=\rho>1$), then one
can construct an irrational number $\alpha$ such that $\{s_n
\alpha\}\in [\delta,1-\delta]$ for all large $n\in\N$, where $\delta>0$
depends on $\rho$ (see \cite{Kat} for
example). As a consequence, the sequence $(s_n)_{n\in\N}$ is not
good for (topological) recurrence. Lastly, we mention that by
considering product systems, one can immediately show that any set
of (topological) recurrence $R$ is partition regular, meaning that
if $R$ is partitioned into finitely many pieces then at least one
of these pieces must still be a set of (topological) recurrence.
Using this observation, one concludes for example that any union
of finitely many lacunary sequences is not a set of recurrence.

We present now  some examples of sets of recurrence:

\begin{theorem} \label{T:examples1} The following are sets of recurrence:
	\begin{enumerate}[(i)]
	\item  Any set of the form
	$\bigcup_{n\in\N}\{a_n,2a_n,\ldots,na_n\}$ where $a_n\in\N$. This follows from a finitary version of Szemer\'edi's theorem.
	
\item Any IP-set, meaning a set that consists of all finite sums
	of some infinite set (\cite{FK12}).
	
\item Any  difference set $S-S$, meaning a set that consists of
	all possible differences of some infinite set $S$.
	
\item The set $\{p(n),n\in \N\}$ where $p$ is any nonconstant
	integer polynomial with $p(0)=0$ (\cite{Fu2, Sa}). In
	fact we only have to assume that the range of the polynomial
	contains multiples of an arbitrary positive integer, this follows from Theorem~\ref{T:single} below.

\item The set $\{p(n),n\in S\}$,  where $p$ is an integer
	polynomial  with $p(0)=0$ and $S$ is any IP-set (\cite{BFM}).
	
\item The set of values of an admissible generalized
	polynomial. This class contains in particular the smallest
	function algebra $G$ containing all integer polynomials
	having zero constant term and such that if $g_1,\ldots ,g_k\in
	G$ and $c_1,\ldots ,c_k\in \R$ then $[\kern-1.5 pt [
	\sum_{i=1}^k c_ig_i]\kern -1.5 pt ]\in G$, where
	$[\kern -1.5 pt [x]\kern -1.5 pt ] =[x+\frac{1}{2}]$ denotes
	the integer nearest to $x$ (\cite{BHaM}).
	
\item The set of shifted primes $\{p-1, p \text{ prime}\}$, and
	the set $\{p+1, p \text{ prime}\}$ $($\cite{Sa}$)$.
	
\item $R=\{ [a(1)], [a(2)],\ldots\} $, where $a(x)=x^c$ for any $c>0$. This follows from Theorem~\ref{T:single} below and standard exponential sum estimates;  see also
	\cite{BKQW05} for a more general result regarding Hardy sequences.
	
\item The set of values of a random non-lacunary
	sequence. More precisely, pick
	$n\in\N$ independently with probability $b_n$ where
	$0\leq b_n\leq 1$ and $\lim_{n\to\infty}nb_n=\infty$, then  the
	resulting set is almost surely a set of recurrence.
	
\item Arbitrary shifts of integers with an even (or odd) number of prime factors and other similar sets with arithmetic structure (\cite{FH15b}).
\end{enumerate}	
\end{theorem}

Showing that the  first three  sets are good for recurrence is  a
straightforward modification  of the argument used to prove
Theorem~\ref{PoincareI}. The other examples
require more work.

A criterion of Kamae and Mend\'es-France (\cite{KM}) provides a
powerful tool that may be used in many instances to establish that
a set $R$ is a set of recurrence. We mention a variation of their
result:
\begin{theorem}[Kamae \& Mend\'es-France~\cite{KM}]\label{T:single}
Suppose that $R=\{a_1<a_2<\ldots\}$ is  a subset of $\N$ such
that:

	\begin{enumerate}[(i)]
		\item  The sequence  $\{a_n \alpha\}_{n\in \N}$ is
uniformly distributed in $\mathbb{T}$ for every irrational
$\alpha$.

\item The set $R_m=\{n\in \N\colon m|a_n\}$ has positive
density for every $m\in \N$.
\end{enumerate}

\noindent Then $R$ is a set of recurrence.
\end{theorem}
We sketch a proof for this result. First, recall Herglotz's
theorem: if $(a_n)_{n\in \Z}$ is a positive definite sequence,
then there is a unique measure $\sigma$ on the torus $\T$ such
that $a_n =\int e^{2\pi i nt}\; d\sigma(t)$. The case of interest
to us is $a_n= \int_{\T} f(x) \cdot  f(T^nx)\; d\mu$, where $T$ is
measure preserving and $f\in L^\infty(\mu)$; $(a_n)$ is positive
definite, and we call $\sigma=\sigma_f$ the {\it spectral measure}
of $f$.

Let now $(X, \mathcal{B},\mu,T)$ be a measure preserving system
and $A\in\mathcal{B}$ with $\mu(A)>0$. Putting $f={\bf 1}_A$, one
has
\begin{equation}\label{123}
\lim_{N\to\infty}\frac{1}{N}\sum_{n=1}^N \int f(x)\cdot
f(T^{a_n}x)\ d\mu=\int_{\T}
\lim_{N\to\infty}\Big(\frac{1}{N}\sum_{n=1}^N e^{2\pi i a_nt}\Big)\
d\sigma_f(t).
\end{equation}
For $t$
irrational the  limit inside the integral is zero (by condition
$(i)$), so the last integral can be taken over the rational points
in $\T$. Since the spectral measure of a function orthogonal to
the subspace
\begin{equation}
\mathcal{H}=\overline{\{f\in  L^2(\mu)\colon \text{ there exists
}k\in \N \text{ with } T^kf=f\}}
\end{equation}
has no
rational point masses, we can easily deduce that when computing
the first limit in \eqref{123}, we can replace the function $f$ by
its orthogonal projection $g$ onto the subspace $\mathcal{H}$ ($g$
is again nonnegative and $g\neq 0$). To complete the argument, we
approximate $g$ by a function $g'$ such that $T^m g'=g'$ for some
appropriately chosen $m$, and use condition $(ii)$ to deduce that
the limit of the average \eqref{123} is positive. 
\medskip

In order to apply Theorem~\ref{T:single}, one uses the standard
machinery of uniform distribution. Recall {\it Weyl's criterion}:
a real valued sequence $(x_n)_{n\in \N}$ is uniformly distributed
mod 1 if for every non-zero $k\in \Z$,
$$
 \lim_{N\to\infty}\frac{1}{N} \sum_{n=1}^N e^{2\pi i k x_n} =0.
$$
This criterion becomes
especially useful when paired with van der
Corput's so-called 
third principal property: if, for every $h\in \N$,
$(x_{n+h}-x_n)_{n\in \N}$ is uniformly distributed mod 1, then
$(x_n)_{n\in \N}$ is uniformly distributed mod $1$. Using the
foregoing criteria and some standard (albeit nontrivial) exponential
sum estimates, one can verify for example that the sets $(iv)$ and
$(vii)$ in Theorem~\ref{T:examples1} are good for recurrence.

In light of the connection elucidated above between uniform
distribution mod 1 and recurrence, it is not surprising that
van der Corput's method has been adapted by modern ergodic
theorists for use in establishing recurrence properties
{\it directly.}

\begin{theorem} [Bergelson~\cite{Be5}]
 Let $(x_n)_{n\in \N}$ be a bounded sequence in a Hilbert space.
 If
 $$\text{D-}\!\!\lim_{m\to\infty}
 \Big( \lim_{N\to\infty}\frac{1}{N} \sum_{n=1}^N \langle
x_{n+m},x_n\rangle\Big) =0,
 $$
then
$$
\lim_{N\to\infty}\norm{\frac{1}{N} \sum_{n=1}^N x_n}=0.
$$
\end{theorem}

Let us illustrate how one uses this ``van der Corput trick'' by
showing that
$S=\{n^2\colon n\in\N\}$ is a set of recurrence. We will actually
establish the following stronger fact: If $(X,\mathcal{B},\mu,T)$
is a measure preserving system and $f\in L^\infty(\mu)$ is
nonnegative and $f\neq 0$  then
\begin{equation}\label{squares}
\liminf_{N\to\infty}\frac{1}{N}\sum_{n=1}^N \int f(x)\cdot
f(T^{n^2}x)\ d\mu>0.
\end{equation}
Then our result follows by setting $f={\bf 1}_A$ for some $A\in
\mathcal{B}$ with $\mu(A)>0$.

The main idea is one that occurs frequently in ergodic theory;
split the function $f$ into two components,  one of which
contributes zero to the limit appearing in \eqref{squares}, and
the other one being much easier to handle than $f$. To do this
consider the $T$-invariant subspace of $L^2(X)$ defined by
\begin{equation}\label{rational}
\mathcal{H}=\overline{\{f\in  L^2(\mu)\colon \text{ there exists
}k\in \N \text{ with } T^kf=f\}}.
\end{equation}
Write $f=g+h$ where $g\in \mathcal{H}$ and $h\bot \mathcal{H}$,
and expand the average in \eqref{squares} into a sum of four
averages involving the functions $g$ and $h$. Two of these
averages vanish because iterates of $g$ are orthogonal to iterates
of $h$. So in order to show that the only contribution comes from
the average that involves the function $g$ alone, it suffices to
establish that
\begin{equation}\label{n^2}
\lim_{N\to\infty}\norm{\
\frac{1}{N}\sum_{n=1}^NT^{n^2}h}_{L^2(\mu)}=0.
\end{equation}
To show this we will apply the Hilbert space van der Corput lemma.
For given $h\in \N$, we let $x_n=T^{n^2}h$ and compute
\begin{align*}
\lim_{N\to\infty}\frac{1}{N} \sum_{n=1}^N \langle
x_{n+m},x_n\rangle&= \lim_{N\to\infty}\frac{1}{N} \sum_{n=1}^N\int
T^{n^2+2nm+m^2}h\cdot T^{n^2}h\ d\mu \\
&=\lim_{N\to\infty}\frac{1}{N} \sum_{n=1}^N\int
T^{2nm}(T^{m^2}h)\cdot h\ d\mu.
\end{align*}
Applying the ergodic theorem to the
transformation $T^{2m}$ and using the fact that $h\bot
\mathcal{H}$, we get that the last limit is $0$. This implies
\eqref{n^2}.

Thus far we have shown that  in order to compute the limit in
\eqref{squares} we can assume that $f=g\in\mathcal{H}$ ($g$ is
also nonnegative and $g\neq 0$). By the definition of
$\mathcal{H}$, given any $\varepsilon>0$, there exists a function
$f'\in \mathcal{H}$ such that $T^kf'=f'$ for some $k\in \N$ and
$\norm{f-f'}_{L^2(\mu)}\leq \varepsilon$. Then
 the limit in \eqref{squares} is  at least $1/k$ times the limit
$$
\liminf_{N\to\infty}\frac{1}{N}\sum_{n=1}^N \int f(x)\cdot
f(T^{(kn)^2}x)\ d\mu.
$$
Applying the triangle inequality twice we get that this is greater
or equal than
\begin{align*}
 \lim_{N\to\infty}\frac{1}{N}\sum_{n=1}^N \int
f'(x)\cdot f'(T^{(kn)^2}x)\ d\mu- c\cdot\varepsilon &=\int
(f'(x))^2
\ d\mu-2\varepsilon \\
&\geq \Big(\int f'(x) \ d\mu\Big)^2-c\cdot \varepsilon,
\end{align*}
for some constant $c$ that does not depend on $\varepsilon$  (we
used that $T^kf'=f'$ and the Cauchy-Schwartz inequality). Choosing
$\varepsilon$ small enough we conclude that the last quantity is
positive, completing the proof.


\section{Multiple Recurrence}\label{S:multiple}
Simultaneous multiple returns of positive measure sets to themselves were
first considered by
H.~Furstenberg~(\cite{Fu1}), who gave a new proof of Szemer\'edi's
theorem~(\cite{Sz}) on arithmetic progressions by deriving it from
the following theorem:
\begin{theorem}[Furstenberg~\cite{Fu1}]\label{T:Furstenberg}
 Let $(X,\mathcal{B},\mu,T)$ be a
measure preserving system and $A\in \mathcal{B}$ with $\mu(A)>0$.
Then for every $k\in\N$, there is some $n\in\N$ such that
\begin{equation}\label{multiple}
 \mu(A\cap
T^{-n}A\cap\cdots\cap T^{-kn}A)>0.
\end{equation}
 \end{theorem}

Furstenberg's proof came by means of a new structure theorem
allowing one to decompose an arbitrary measure preserving system
into component elements exhibiting one of two extreme types of
behavior: {\it compactness}, characterized by regular, ``almost
periodic'' trajectories, and {\it weak mixing}, characterized by
irregular, ``quasi-random'' trajectories. On $\T$, these
types of behavior are exemplified by rotations and by the doubling
map, respectively. To see the point, imagine trying to predict the
initial digit of the dyadic expansion of $T^nx$ given knowledge of
the initial digits of $T^ix$, $1\leq i <n$. We use the case $k=2$
to illustrate the basic idea.

It suffices to
show that if $f\in L^\infty(\mu)$ is nonnegative and $f\neq 0$,
one has
\begin{equation}\label{Roth}
\liminf_{N\to\infty}\frac{1}{N}\sum_{n=1}^N \int f(x)\cdot
f(T^nx)\cdot f(T^{2n}x)\ d\mu>0.
\end{equation}
An ergodic decomposition argument
enables us to assume that our system is ergodic.
As in the earlier case of the squares, we split $f$ into ``almost
periodic'' and ``quasi-random'' components. Let $\mathcal{K}$ be the
closure in $L^2$ of the subspace spanned by the eigenfunctions of
$T$, i.e. the functions $f\in L^2(\mu)$ that satisfy
$f(Tx)=e^{2\pi i \alpha} f(x)$ for some $\alpha\in \R$. We write
$f=g+h$, where $g\in \mathcal{K}$ and $h\bot \mathcal{K}$. It can
be shown that $g,h\in L^{\infty}(\mu)$ and $g$ is again
nonnegative with $g\neq 0$. We expand the average in \eqref{Roth}
into a sum of eight averages involving the functions $g$ and $h$. In
order to show that the only non-zero contribution to the limit
comes from the term involving $g$ alone, it suffices to establish
that
\begin{equation}\label{n,2n}
\lim_{N\to\infty}\norm{\ \frac{1}{N}\sum_{n=1}^N T^ng \cdot T^{2n}h
}_{L^2(\mu)}=0,
\end{equation}
(and similarly with $h$ and $g$ interchanged, and with $g=h$,
which is similar). To establish \eqref{n,2n}, we use the Hilbert
space van der Corput lemma on $x_n=T^{n}g\cdot T^{2n}h$. Some
routine computations and a use of the ergodic theorem reduce the
task to showing that
$$
\text{D-}\!\!\lim_{m\to\infty}\Big(\int h(x)\cdot h(T^{2m}x)\
d\mu\Big)=0.
$$
But this is well known for $h\bot \mathcal{K}$ (in
virtue of the fact that
for $h\bot \mathcal{K}$ the spectral measure $\sigma_h$ is
continuous, for example).

We are left with the average \eqref{Roth} when
$f=g\in \mathcal{K}$. In this case $f$ can be approximated
arbitrarily well by a linear combination of eigenfunctions,
which easily implies that given $\varepsilon>0$ one has $\norm
{T^nf-f}_{L^2(\mu)}\leq \varepsilon$ for a set of $n\in\N$ with
bounded gaps. Using this fact and the triangle inequality,
one finds that for a set of $n\in\N$ with bounded gaps,
$$
 \int f(x)\cdot
f(T^nx)\cdot f(T^{2n}x)\ d\mu\geq \Big(\int f\ d\mu\Big)^3 -c\cdot
\varepsilon
$$
for a constant $c$ that is independent of $\varepsilon$. Choosing
$\varepsilon$ small enough, we get \eqref{Roth}.
\medskip

The new techniques developed for the proof of Theorem
\ref{T:Furstenberg} have led to a number of extensions, many of which
have to date only ergodic proofs.
To expedite discussion of some of these
developments, we introduce a definition:
 \begin{definition}
Let   $R\subset \Z$ and $k\in\N$. Then $R$ is a set of $k$-{\em
recurrence}  if for every invertible measure preserving system
$(X,\mathcal{B},\mu,T)$ and $A\in \mathcal{B}$ with $\mu(A)>0$,
there is some nonzero $n\in R$ such that
$$
\mu(A\cap T^{-n}A\cap\cdots\cap T^{-kn}A)>0.
$$
\end{definition}
The notions of $k$-recurrence are distinct for different values of $k$.
An example of a difference set that is a set of
$1$-recurrence but not a set of $2$-recurrence was given in \cite{Fu1};
sets of $k$-recurrence that are not sets of $(k+1)$-recurrence for general
$k$ were given in \cite{FLW} ($R_k=\{n\in\N\colon \{n^{k+1}\sqrt{2}\}\in
[1/4,3/4]\}$ is such).

Aside from difference sets, the sets of ($1$-)recurrence given in
Theorem~\ref{T:examples1} may well
be  sets of $k$-recurrence for every $k\in \N$, though this has not
been verified in all cases. Let us summarize
the current state of knowledge.
The following are sets of $k$-recurrence for every $k$:
Sets of the form $\bigcup_{n\in\N}\{a_n,2a_n,\ldots,na_n\}$
where $a_n\in\N$ (this follows from a uniform version of
Theorem~\ref{T:Furstenberg} that can be found in \cite{BHRF}).
Every IP-set (\cite{FK12}). The set $\{p(n),n\in \N\}$ where $p$
is any nonconstant integer polynomial with $p(0)=0$ (\cite{BL}),
and more generally, when the range of the polynomial contains
multiples of an arbitrary integer (\cite{Fr}). The set
$\{p(n),n\in S\}$ where $p$ is an integer polynomial with $p(0)=0$
and $S$ is any IP-set (\cite{BM}). The set of values of an
admissible generalized polynomial (\cite{BM10,M}). Moreover,   it was shown in \cite{FHK} for $k=2$ and in \cite{WZ11} for general $k\in\N$,
that  the set
of shifted primes $\{p-1, p \text{ prime}\}$, and  the set $\{p+1,
p \text{ prime}\}$ are sets of $k$-recurrence (see also \cite{BLZ11,FrHK11,Kou15a,Su15} for related work). Several other  multiple recurrence
results were obtain in the last ten years,  including results for Hardy sequences (\cite{BMR,Fr09,Fr10,FrW09}), various subsets of integer part polynomial sequences  (\cite{KK,Kou15a,Kou15b}), random sequences (\cite{FrLW11,FrLW14}), and sets of arithmetic nature (\cite{BKL,FH15a,FH15b}).

 More generally, one would like to know for which sequences of
integers  $a_1(n)$,$\ldots$,$a_k(n)$ it is the case that for every
invertible measure preserving system $(X,\mathcal{B},\mu,T)$ and
$A\in \mathcal{B}$ with $\mu(A)>0$, there is some nonzero $n\in
\N$ such that
\begin{equation}
\mu(A\cap T^{-a_1(n)}A\cap\cdots\cap T^{-a_k(n)}A)>0.
\end{equation}
Unfortunately, a criterion analogous to the one given in
Theorem~\ref{T:single} for $1$-recurrence is not yet available for
$k$-recurrence when $k>1$. Nevertheless, there have been some
notable positive results, such as the following:
\begin{theorem}[Bergelson \& Leibman~\cite{BL}]\label{T:BL}
 Let $(X,\mathcal{B},\mu,T)$ be an invertible
measure preserving system and $p_1(n),\ldots,p_k(n)$ be integer
polynomials with zero constant term. Then for every $A\in
\mathcal{B}$ with $\mu(A)>0$, there is some $n\in\N$ such that
\begin{equation}\label{multiple2}
 \mu(A\cap
T^{-p_1(n)}A\cap\cdots\cap T^{-p_k(n)}A)>0.
\end{equation}
 \end{theorem}
Furthermore, it has been shown that the $n$ in \eqref{multiple2}
can be chosen from any IP set (\cite{BM}), and the polynomials
$p_1,\ldots, p_k$ can be chosen to belong to the more general
class of admissible generalized  polynomials~(\cite{M}) or the class of intersective polynomials~(\cite{BLL}).

An important boost in the area of multiple recurrence was
given by a breakthrough of Host and Kra (\cite{HK1}).
Building on work of Conze and Lesigne (\cite{CL1, CL2}) and
Furstenberg and Weiss (\cite{FW}) (see also the excellent
survey \cite{Kra}, exploring close parallels with \cite{GT1}
and the seminal paper of
Gowers~(\cite{Gow})), they
isolated the structured component (or factor) of a measure
preserving system that one needs to analyze in order to prove
several multiple recurrence and convergence results. This allowed
them, in particular, to prove existence of $L^2$ limits for the
so-called ``Furstenberg ergodic averages'' ${1\over N}\sum_{n=1}^N
\prod_{i=0}^k f(T^{in}x)$, which had been a major open problem
since the original ergodic proof of Szemer\'edi's theorem.
Subsequently Ziegler in \cite{Zi}  gave a new proof of the
aforementioned limit theorem and established minimality of the
factor in question. It turns out that this minimal component
admits of a purely algebraic characterization; it is a {\em
nilsystem}, i.e. a  rotation on a homogeneous space of a nilpotent
Lie group. This fact, coupled with some recent results about
nilsystems (see \cite{Le1, Le2} for example), makes the
analysis of some otherwise intractable multiple recurrence
problems much more manageable. These developments
have made it possible to obtain new multiple recurrence results and they also allowed us to
estimate the size of the multiple
intersection in \eqref{multiple} for $k=2,3$ (the case $k=1$ is
Theorem~\ref{Khin}):

\begin{theorem}[Bergelson, Host \& Kra
\cite{BHK}]\label{T:BHK}
Let $(X,\mathcal{B},\mu,T)$ be an ergodic
measure preserving system and $A\in \mathcal{B}$. Then for
$k=2,3$ and for every $\varepsilon>0$,
\begin{equation}\label{bounds}
\mu(A\cap T^{-n}A\cap\cdots\cap T^{-kn}A)>\mu^{k+1}(A)-\varepsilon
\end{equation}
for a set of $n\in\N$ with bounded gaps.
\end{theorem}
Based on work of Ruzsa that appears as an appendix to the paper,
it is also shown in \cite{BHK} that a similar estimate fails for
ergodic systems (with any power of $\mu(A)$ on the right hand
side) when $k\geq 4$. Moreover, when the system is nonergodic it
also fails for $k=2,3$, as can be seen with the help of an example
in \cite{Beh}. Again considering the doubling map $Tx=2x$
and the set $A=[0,1/2]$,
one sees that the positive results for $k\leq 3$ are sharp.  When
the polynomials $n,2n,\ldots,kn$ are replaced by linearly
independent polynomials $p_1,p_2,\ldots,p_k$ with zero constant
term, similar lower bounds hold for every $k\in \N$ without
assuming ergodicity (\cite{FK2}). The case where the polynomials $n,2n,3n$ are
replaced with general polynomials $p_1,p_2,p_3$ with zero constant
term is treated in \cite{Fr} (see also \cite{DLMS}) and more general results involving Hardy field sequences and polynomials evaluated at the primes
are obtained in \cite{DLMS}.

\section{Connections with Combinatorics and Number
Theory}\label{applications} The combinatorial ramifications of
ergodic-theoretic
recurrence were first observed by Furstenberg, who perceived a
correspondence between recurrence properties of measure preserving
systems and the existence of structures in sets of integers having
positive upper density. This gave rise to the field of ergodic
Ramsey theory, in which problems in combinatorial number theory
are treated using techniques from ergodic theory. The
following formulation is from \cite{Be527}.

\begin{theorem}\label{correspondence} Let $\Lambda$ be a
subset of the integers. There exists an invertible measure
preserving system $(X,\B,\mu,T)$ and a set $A\in\mathcal{B}$ with
$\mu(A)=\overline{d}(\Lambda)$ such that
\begin{equation}\label{E:correspondence}
\overline{d}(\Lambda \cap (\Lambda-n_1)\cap\ldots\cap
(\Lambda-n_k))\geq \mu(A\cap T^{-n_1}A\cap\cdots \cap T^{-n_k}A),
\end{equation}
for all $k\in\mathbb{N}$ and  $n_1,\ldots,n_k\in\Z$.
\end{theorem}
\begin{proof}
The space $X$ will be taken to be the sequence space $\{0,1\}^\Z$,
$\mathcal{B}$ is the Borel $\sigma$-algebra,  while $T$ is the
shift map defined by $(Tx)(n)=x(n+1)$ for $x\in \{0,1\}^\Z$, and
$A$ is the set of sequences $x$ with $x(0)=1$. So the only thing
that depends on $\Lambda$ is the measure $\mu$ which we now
define. For $m\in\N$ set $\Lambda^0=\Z\setminus \Lambda$ and
$\Lambda^1=\Lambda$. Using a
diagonal argument we can find an increasing sequence of integers
$(N_m)_{m\in\N}$ such that $\lim_{m\to\infty}
|\Lambda\cap [1,N_m]|/N_m=\overline{d}(\Lambda)$ and such that
\begin{equation}\label{E:limit}
\lim_{m\to\infty} \frac{|(\Lambda^{i_1}-n_1)\cap
(\Lambda^{i_2}-n_2)\cap\cdots\cap (\Lambda^{i_r}-n_r)\cap [1,N_m]
|}{N_m}
\end{equation}
exists for every $n_1,\ldots, n_r\in\bbZ$, and
$i_1,\ldots,i_r\in\{0,1\}$. For $n_1, n_2, \ldots, n_r\in\bbZ$, and $i_{1}, i_2, \ldots,
i_{r}\in\{0,1\}$, we define the measure $\mu$ of the {\it cylinder set}
$\{x_{n_1}=i_{1},x_{n_2}=i_{2},\ldots,x_{n_r}=i_{r}\}$ to be the
limit \eqref{E:limit}.
Thus defined, $\mu$ extends to a premeasure on the algebra of sets
generated by cylinder sets and hence by Carath\'eodory's extension
theorem (\cite{Ca})  to a probability measure on $\mathcal{B}$. It
is easy to check that $\mu(A)=\overline{d}(\Lambda)$, the shift
transformation $T$ preserves the measure $\mu$ and
\eqref{E:correspondence} holds.
\end{proof}

Using this principle for $k=1$, one may check that any set of
recurrence is {\it intersective}, that is intersects $E-E$ for
every set $E$ of positive density.
Using it for
$n_1=n,n_2=2n,\ldots, n_k=kn$, together with
Theorem~\ref{T:Furstenberg}, one gets an ergodic proof of
Szemer\'edi's theorem~(\cite{Sz}), stating that every subset of the integers
with positive  upper density contains arbitrarily long arithmetic
progressions (conversely, one can easily deduce
Theorem~\ref{T:Furstenberg} from  Szemer\'edi's theorem, and that
intersective sets are sets of recurrence). Making the choice
$n_1=n^2$ and using part $(iv)$ of Theorem~\ref{T:single}, we get
an ergodic proof of   the surprising result of S\'ark\"ozy
(\cite{Sa}) stating  that every subset of the integers with
positive upper density contains two elements whose difference is a
perfect square. More generally, using
Theorem~\ref{correspondence}, one can translate all of the
recurrence results of the previous two sections to results in
combinatorics. (This is not straightforward for
Theorem~\ref{T:BHK} because of the ergodicity assumption made
there. We refer the reader to \cite{BHK} for the combinatorial
consequence of this result.) We mention explicitly only the
combinatorial consequence of Theorem~\ref{T:BL}:
\begin{theorem}[Bergelson \& Leibman ~\cite{BL}]
Let $\Lambda\subset\Z$ with $\overline{d}(\Lambda)>0$, and
$p_1,\ldots,p_k$ be integer polynomials with zero constant term.
Then $\Lambda$ contains infinitely many configurations of the form
$\{x,x+p_1(n),\ldots,x+p_k(n)\}$.
\end{theorem}
The ergodic  proof is the only one known for this
result, even for patterns of the form $\{x,x+n^2,x+2n^2\}$ or
$\{x,x+n,x+n^2\}$.

Ergodic-theoretic contributions to the field of geometric Ramsey
theory were made by Furstenberg, Katznelson, and Weiss (\cite{FKW}), who
showed that if $E$ is a positive upper density
subset of $\R^2$ then: $(i)$  $E$ contains
points with any large enough distance (see also \cite{Bo0} and
\cite{FM}), $(ii)$ Every $\delta$-neighborhood of $E$ contains
three points forming a triangle congruent to any given large
enough dilation of a given triangle (in \cite{Bo0} it is shown
that if the three points lie on a straight line one cannot always
find three points with this property in $E$ itself). Recently, a
generalization of property $(ii)$  to arbitrary finite
configurations of $\R^m$ was obtained by Ziegler~(\cite{Zi2}).

We  also  mention some  exciting connections of
multiple recurrence  with some structural properties of the set of
prime numbers. The first one is in the work of Green and
Tao~(\cite{GT1}), where the existence of arbitrarily long
arithmetic progressions of primes was demonstrated, the authors,
in addition to using Szemer\'edi's theorem outright, use several
ideas from its ergodic-theoretic proofs, as appearing in
\cite{Fu1} and \cite{FKO}. The second one is in the recent work of
Tao and Ziegler~\cite{TZ}, a quantitative version of
Theorem~\ref{T:BL} was used to prove that the primes contain
arbitrarily long polynomial progressions. Furthermore,
 results in ergodic theory related to the structure of the
minimal characteristic factors of certain multiple ergodic
averages, play an important role in the work of Green,  Tao, and Ziegler
(\cite{GT2,GT3,GT4,GT09d, GTZ}), were they   get asymptotic formulas for the number of $k$-term
arithmetic progressions of primes up to $x$. This work verifies
an interesting special case of the
Hardy-Littlewood $k$-tuple conjecture predicting the asymptotic
growth rate of $N_{a_1,\ldots ,a_k}(x)=$ the number of
configurations of primes having the form
$\{p,p+a_1,\ldots,p+a_k\}$ with $p\leq x$.

In a more recent development, the tools developed in the last two decades to deal with delicate multiple recurrence problems
have played an instrumental role in analyzing the structure of measure-preserving
systems naturally associated with bounded multiplicative functions. These results were used in the last two
years in works of Tao and Ter\"av\"ainen (\cite{TT18,TT19}) to make progress on the Chowla and Elliott conjectures and
in works of Frantzikinakis and Host (\cite{FH18,FH19})  to make progress on the M\"obius disjointness conjecture of Sarnak.
It appears that  this interplay of ergodic theory and number theory is going to be  essential for the resolution of
 several notoriously difficult problems
concerning higher order correlations of multiplicative and other number theoretic functions.

Finally, we remark that in this article we have restricted
attention to multiple recurrence and Furstenberg correspondence
for $\Z$ actions, while in fact there is a wealth of literature on
extensions of these results to general commutative, amenable and
even non-amenable groups. For an excellent exposition of these and
other recent developments the reader is referred to the survey articles
\cite{Au16, Be1, Be2, Be06b}. Here, we give just one notable combinatorial
corollary to some work of this kind, a density version of the
classical Hales-Jewett coloring theorem (\cite{HJew}).

\begin{theorem}[Furstenberg \&
Katznelson ~\cite{FKat2}; see also \cite{Poly1}]\label{T:FK} Let $W_n(A)$ denote the set
of words of length $n$ with letters in the alphabet
$A=\{a_1,\ldots,a_k\}$. For every $\varepsilon>0$ there exists
$N_0=N_0(\varepsilon,k)$ such that if $n\geq N_0$ then any subset
$S$ of $W_n(A)$ with $|S|\geq \varepsilon k^n$ contains a
combinatorial line, i.e., a set consisting of $k$ $n$-letter words,
having fixed letters in $l$ positions, for some $0\leq l<n$, the remaining
$n-l$ positions being occupied by a variable letter $x$, for
$x=a_1,\ldots,a_k$. (For example, in $W_4(A)$ the sets $\{
(a_1,x,a_2,x)\colon x\in A\}$ and $\{(x,x,x,x),\colon x\in A\}$
are combinatorial lines.)
\end{theorem}
At first glance, the uninitiated
reader may not appreciate the importance of this ``master'' density
result, so it is instructive to derive at least one of its
immediate consequences. Let $A=\{0,1,\ldots,k-1\}$ and interpret
$W_n(A)$ as integers in base $k$ having at most $n$ digits. Then a
combinatorial line in $W_n(A)$ is an arithmetic progression of
length $k$--for example, the line $\{(a_1,x,a_2,x)\colon
x\in A\}$ corresponds to the progression $\{ m, m+n, m+2n, m+3n\}$,
where $m=a_1+a_2d^2$ and $n=d+d^3$. This allows one to deduce
Szemer\'edi's theorem. Similarly, one can deduce from
Theorem~\ref{T:FK} multidimensional and IP extensions of
Szemer\'edi's theorem~(\cite{FK1, FK12}), and some related
results about vector spaces over finite fields~(\cite{FK12}).

\section{Future Directions}\label{problems}
In this section we formulate a few open problems relating to the
material in the previous three sections. It should be noted that
this selection reflects the authors' interests, and does not strive
for completeness. A more extensive list of problems related to ergodic theory of $Z$-actions can be found in \cite{Fr16}.

We start with an intriguing question of Katznelson~(\cite{Kat})
about sets of topological recurrence. A set $S\subset\N$ is a set
of Bohr recurrence if for every $\alpha_1,\ldots,\alpha_k\in \R$
and $\varepsilon>0$ there exists $s\in S$ such that
$\{s\alpha_i\}\in [0,\varepsilon]\cup [1-\varepsilon,1)$ for
$i=1,\ldots,k$.
\begin{problem}
Is every set of Bohr recurrence a set of topological recurrence?
\end{problem}
Background for this problem and evidence for a positive answer can
be found in \cite{Kat, W}. A negative for a related question concerning
general Abelian groups is given in \cite{Gr16}. As we mentioned in
Section~\ref{single}, there exists a set of topological recurrence
(and hence Bohr recurrence) that is not a set of recurrence.
\begin{problem}
Is the set $S=\{l!2^m3^n\colon l,m,n\in \N\}$ a set of recurrence?
Is it a set of $k$-recurrence for every $k\in \N$?
\end{problem}
It can be shown that $S$ is a set of Bohr recurrence.
Theorem~\ref{T:single} cannot be applied since the uniform
distribution condition fails for some irrational numbers $\alpha$.
A relevant question was asked by Bergelson in \cite{Be1}: ``Is the
set $S=\{2^m3^n\colon m,n\in\N\}$ good for single recurrence for
weakly mixing systems?''

We mentioned in Section~\ref{single} that the sets of shifted primes and the set of fractional powers of integers are known to be sets of $k$-recurrence for every $k\in \N$. The following related question remains  open:
\begin{problem}
	Show that for every $c\in \R^+\setminus \Z$ the set $\{[p^{c}], p\in \mathbb{P}\}$, where $\mathbb{P}$ is the set of primes, is a set of $k$-recurrence for
	every $k\in \N$.
	\end{problem}
The problem is open  even for $k=2$ and $c=\frac{3}{2}$.


We mentioned in Section~\ref{single} that random non-lacunary
sequences (see definition there) are almost surely sets of recurrence.
\begin{problem}
Show that random non-lacunary sequences are almost surely sets of
$k$-recurrence for every $k\in \N$.
\end{problem}
Some progress is made in  \cite{BGSZ16, BDG16, BG17, Chr11, FrLW11, FrLW14} but
the problem is open  even for $k=2$ for random sequences with at most quadratic growth.
 We refer the
reader to the survey  \cite{RW} for a nice exposition of the
argument used by  Bourgain~\cite{Bo} to handle the case $k=1$.

It was shown in  \cite{FLW} that if $S$ is a set of $2$-recurrence
then the set of its squares  is a set of recurrence for circle
rotations. The same method shows that it is actually a set of Bohr
recurrence.
\begin{problem}
If $S\subset \Z$ is a set of $2$-recurrence, is it true that
$S^2=\{s^2\colon s\in S\}$ is a set of  recurrence?
\end{problem}
A similar question was asked in \cite{BGL}: ``If S is a set of
$k$-recurrence for every $k$, is the same true of $S^2$?''.

One would like to find a criterion that would allow one to deduce
that a sequence is good for double (or higher order) recurrence
from some uniform distribution properties of this sequence.
\begin{problem}
Find  necessary conditions for double recurrence similar to the
one given in Theorem~\ref{T:single}.
\end{problem}
It is now well understood that such a criterion should involve
uniform distribution properties of some generalized polynomials or
$2$-step nilsequences.

We mentioned in Section~\ref{applications} that
every positive density subset of $\R^2$ contains points with any
large enough distance. Bourgain (\cite{Bo0}) constructed a
positive density subset $E$  of $\R^2$, a triangle $T$, and numbers $t_n\to\infty$,
such that $E$ does not contain congruent copies of all $t_n$-dilations of $T$.
But the triangle $T$ used in this construction is degenerate,
which leaves the following question open:
\begin{problem}
Is it true that every positive density subset of $\R^2$ contains a
triangle congruent to any large enough dilation of a given
non-degenerate triangle?
\end{problem}
For further discussion on this question the reader can consult the
survey \cite{Gr}.

The following question of
Aaronson and Nakada (\cite{Aa0}) is related to a classical
question of Erd\H{o}s concerning whether every
$K\subset \N$ such that $\sum_{n\in K}1/n=\infty$ contains
arbitrarily long arithmetic progressions:
\begin{problem}
Suppose that $(X,\mathcal{B},\mu,T)$ is a
$\{1/n\}$-conservative ergodic measure preserving system. Is it
true that for every $A\in \mathcal{B}$ with $\mu(A)>0$ and $k\in
\N$  we have $\mu(A\cap T^{-n}A\cap\cdots\cap T^{-kn}A)>0$ for
some $n\in \N$?
\end{problem}
The answer is positive for the class of Markov shifts, and it is
remarked in \cite{Aa0} that  if the Erd\H{o}s conjecture is true
then the answer will be positive in general. The converse is not
known to be true. For a  related result showing that multiple
recurrence is preserved by extensions of infinite measure
preserving systems see \cite{Me}.

Our next problem is motivated by the question whether
Theorem~\ref{T:FK} has a polynomial version (for a precise
formulation of the general conjecture see \cite{Be1}). Not even
this most special consequence of it is known to hold.
\begin{problem}
Let $\varepsilon>0$. Does there exist $N=N(\varepsilon)$ having
the property that every family $P$ of subsets of $\{1,\ldots
,N\}^2$ satisfying $|P|\geq \varepsilon 2^{N^2}$ contains a
configuration $\{ A, A\cup (\gamma\times \gamma)\}$, where
$A\subset \{1,\ldots ,N\}^2$ and $\gamma\subset \{1,\ldots ,N\}$
with $A\cap (\gamma\times \gamma)=\emptyset$?
\end{problem}

A measure preserving action of a general countably infinite group
$G$ is a function $g\rightarrow T_g$ from $G$ into the space of
measure preserving transformations of a probability space $X$ such
that $T_{gh}=T_gT_h$. It is easy to show that a version of
Khintchine's recurrence theorem holds for such actions: if
$\mu(A)>0$ and $\varepsilon >0$ then $\{ g:\mu(A\cap T_g
A)>\big( \mu(A)\big) ^2-\varepsilon\}$ is syndetic.
However it is unknown whether
the following ergodic version of Roth's theorem holds.

\begin{problem}
Let $(T_g)$ and $(S_g)$ be measure preserving $G$-actions of a
probability space $X$ that commute in the sense that
$T_gS_h=S_hT_g$ for all $g,h\in G$. Is it true that for all
positive measure sets $A$, the set of $g$ such that $\mu(A\cap T_g
A\cap S_g A)>0$ is syndetic?
\end{problem}

We remark that for
general (possibly amenable)
groups $G$ not containing arbitrarily large finite subgroups
nor elements of infinite order, it is not known whether one can
find a {\it single} such $g\neq e$.
On the other hand, the answer is known to be positive for general $G$ in case
$(T_{g}^{-1} S_g)$ is a $G$-action (\cite {BM5}); even under such strictures,
however, it is unknown whether a triple recurrence
theorem holds.


\end{document}